\documentclass[submission]{FPSAC2021}

\newtheorem{thm}{Theorem}

\newtheorem{definition}{Definition}
\newtheorem{ex}{Example}

\usepackage{lipsum}
\usepackage{mathtools}
\usepackage{float}
\usepackage{makecell}

\newcommand{\R}{\mathbb{R}}
\newcommand{\Z}{\mathbb{Z}}
\newcommand{\N}{\mathbb{N}}

\title[Double boxes and double dimers]{Double boxes and double dimers}

\author[Tatyana Benko and Benjamin Young]{\texorpdfstring{Tatyana Benko\addressmark{1}\thanks{tbenko@uoregon.edu. Tatyana Benko was partially funded by NSF grant DMS-2039316.} \and Benjamin Young\addressmark{1}\thanks{bjy@uoregon.edu}}{tbenko@uoregon.edu. Tatyana Benko was partially funded by NSF grand DMS-2039316. bjy@uoregon.edu}}

\address{\addressmark{1}Department of Mathematics, University of Oregon, Eugene, OR, USA \\}

\received{\today}

\abstract{We give a combinatorial proof of a result in rank 2 Donaldson-Thomas theory, which states that the generating function for certain plane-partition-like objects, called double-box configurations, is equal to a product of MacMahon's generating function for (boxed) plane partitions. In our proof, we first give the correspondence between double-box configurations and double-dimer configurations on the hexagon lattice with a particular tripartite node pairing. Using this correspondence, we apply graphical condensation and double-dimer condensation to prove the result.}

\keywords{Plane partitions, double-box configurations, dimer model, Donaldson-Thomas theory, condensation}

\begin{document}

\maketitle

\section{Introduction}

In this extended abstract, we enumerate certain plane-partition-like objects called \emph{double-box configurations}.  These objects were introduced by Gholampour, Kool and Young \cite{gholampour2018rank} for the purpose of computing the rank 2 Donaldson-Thomas (DT) invariants of a Calabi-Yau threefold. We define double-box configurations in Definition \ref{def:double_box_config}, as well as their generating function, $Z_{DBC}^{a,b,c}(q)$, in Equation \ref{eq:DBC_GF}, and we give a combinatorial proof of the following geometrically motivated theorem:

\begin{thm}\label{thm:main}
Let $a,b,c \in \N$, then

    \begin{equation}\label{eq:main}
        Z_{DBC}^{a,b,c}(q) = M(q)^2M_{a,b,c}(q)
    \end{equation}
    
    \noindent where $Z_{DBC}^{a,b,c}(q)$ denotes the generating function for double-box configurations, and

    \begin{equation*}
        M(q) = \prod_{i=1}^{\infty} \frac{1}{(1-q^i)^i}
    \end{equation*}

    \noindent is MacMahon's generating function for plane partitions, and
    
    \begin{equation*}
        M_{a,b,c}(q) = \prod_{s=1}^a \prod_{t=1}^b \prod_{r=1}^c \frac{1-q^{s+t+r-1}}{1-q^{s+t+r-2}}
    \end{equation*}
   
    \noindent is MacMahon's generating function for boxed $a \times b \times c$ plane partitions. 
\end{thm}

Note that this formula already has a geometric proof \cite{gholampourkool}. In this paper, we outline our combinatorial proof of this formula, which surprisingly uses the tripartite double-dimer model of Kenyon and Wilson \cite{kenyonwilson}. We follow the general strategy of \cite{jwy}, where the tripartite double-dimer model also appears (though for apparently completely different reasons). Our proof consists of two main components. First, we give a correspondence between double-box configurations and tripartite double-dimer configurations on the hexagon lattice; this correspondence is many-to-one, yet still weight preserving. Next, we use a quadratic recurrence called condensation relation to prove our main result. We show that to the left hand side of Equation \ref{eq:main}, $Z_{DBC}^{a,b,c}(q)$, we may apply a result by Jenne \cite{jenne2021combinatorics}, which states that under certain conditions the generating function for tripartite double-dimer configurations satisfies a recurrence relation related to the Desnanot-Jacobi identity from linear algebra. Then, using Kuo condensation \cite{kuo2004applications} (also related to the Desnanot-Jacobi identity), we show that $M(q)^2M_{a,b,c}(q)$ satisfies the same recurrence relation. Finally, we show that both sides of Equation \ref{eq:main} satisfy the same initial conditions. 

The full version of this abstract will appear in a longer paper; proofs and some details have been omitted here due to space constraints.

\section{Definitions}

\begin{figure}[H]
\centering
\includegraphics[width=\textwidth]{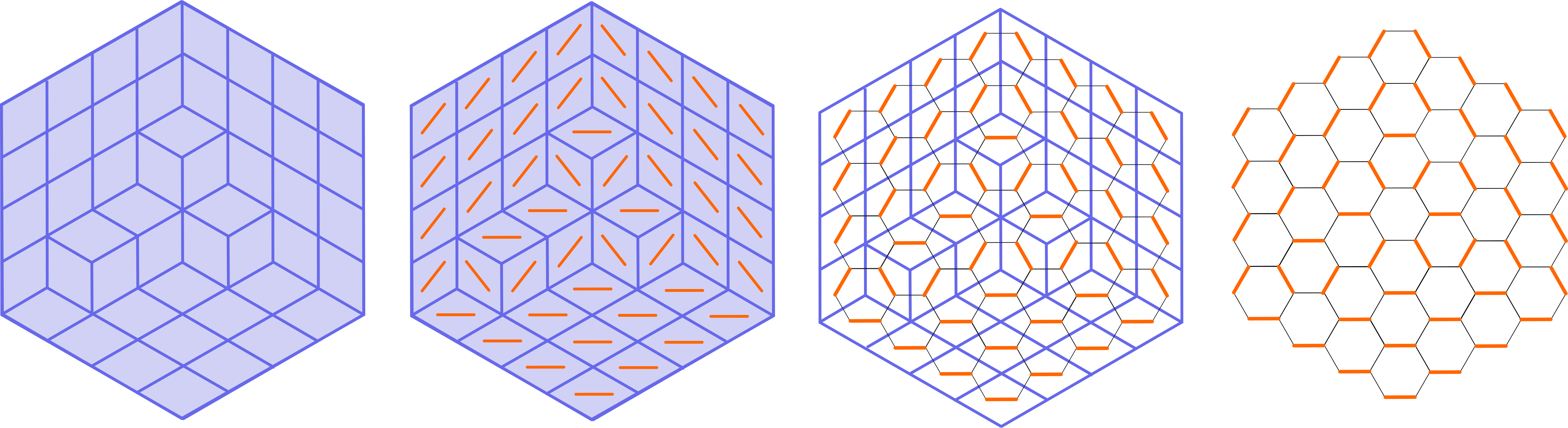}
\caption{Folklore bijection between plane partitions (left-most) and dimer configurations on the hexagon graph (right-most)}
\label{fig:folklore_bij}
\end{figure}

A \emph{plane partition} is a two-dimensional array of nonnegative integers $\pi_{i,j}$ for $i,j \geq 0$, with $\pi_{i,j+1} \geq \pi_{i,j}$ and $\pi_{i+1, j} \geq \pi_{i,j}$ for all $i,j$, with finitely many $\pi_{i,j}$ being nonzero. We can visualize a plane partition $\pi$ as a stack of boxes in the corner of a room, with the number of boxes in each stack given by the entries of $\pi$. A \emph{dimer configuration} (also called a \emph{perfect matching}) on a graph $G = (V,E)$ is a collection of edges $E' \subseteq E$ such that every vertex in $V$ is covered exactly once. There is a bijection between plane partitions and dimer configurations on the hexagon graph, sometimes referred to as the \emph{folklore bijection}  (see Figure \ref{fig:folklore_bij}). The stacks of boxes representing a plane partition $\pi$ can be viewed as a lozenge tiling of a hexagonal region of triangles, corresponding to the visible faces of the boxes and the tiles on the walls and floor of the room. Note that the triangular lattice is dual to the hexagon lattice. Moreover, each lozenge is made of two equilateral triangles that share an edge. If we join the centers of these two triangles with the corresponding dual edge, and do this for all tiles in the tiling, we get a perfect matching on the hexagon graph. 

\begin{figure}[H]
\centering
\includegraphics[width=\textwidth]{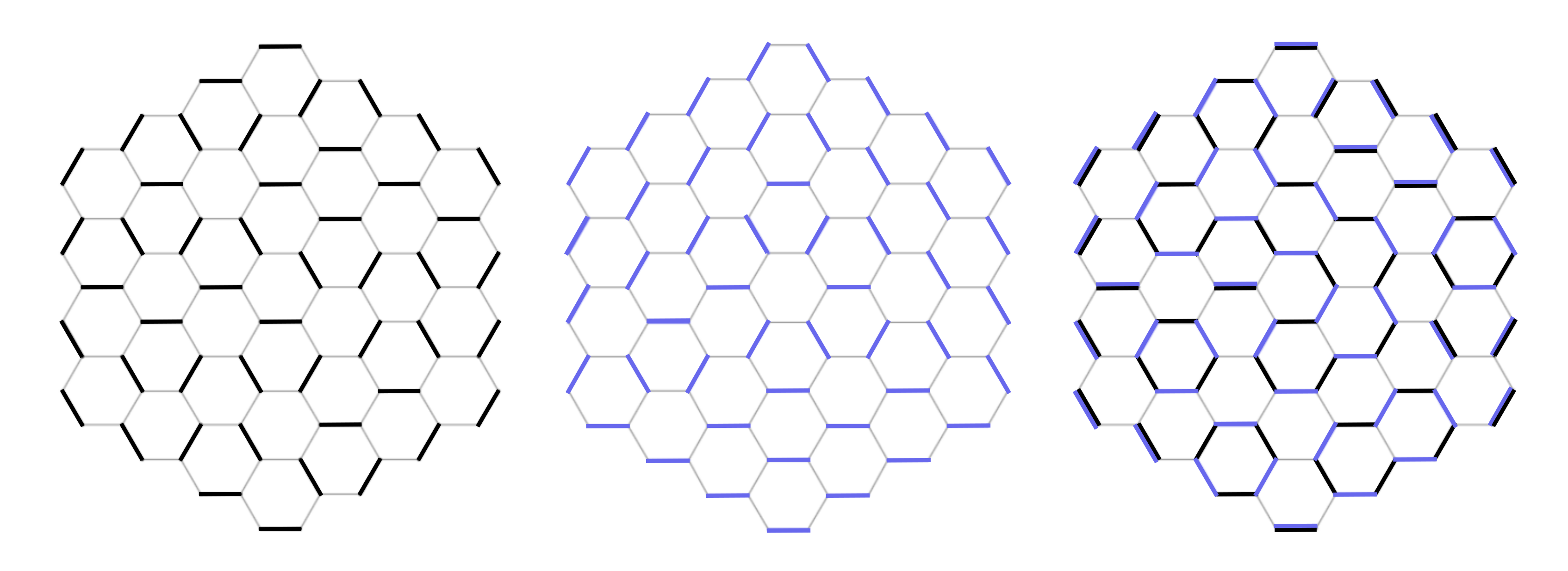}
\caption{Double-dimer configuration on the hexagon graph (right-most) from two single-dimer configurations}
\label{fig:double_dimer}
\end{figure}

Overlaying two perfect matchings of a graph $G=(V,E)$ gives a \emph{double-dimer configuration}, which consists of doubled edges and loops (see Figure \ref{fig:double_dimer}). If in addition we have defined a set of nodes $N \subset V$, that is, a special set of vertices, then the double-dimer configuration on $G=(V,E)$ with node set $N$ is a collection of edges $E' \subseteq E$ such that each vertex in $V \setminus N$ is covered exactly twice, and each node in $N$ is covered exactly once. In this case, the double-dimer configuration consists of doubled edges, loops, and paths between the nodes in $N$. 

\subsection{Double-box configurations}

In this section we define the double-box configurations, and provide some examples. Throughout this section we let $a,b,c \in \mathbb{N}$ be fixed. We identify the point $(i,j,k) \in \Z^3$ with the unit cube

$$[i,i+1] \times [j,j+1] \times [k,k+1] \in \R^3.$$

\noindent We refer to this unit cube as the box $(i,j,k)$.

\begin{figure}[H]
\centering
\includegraphics[width=0.7\textwidth]{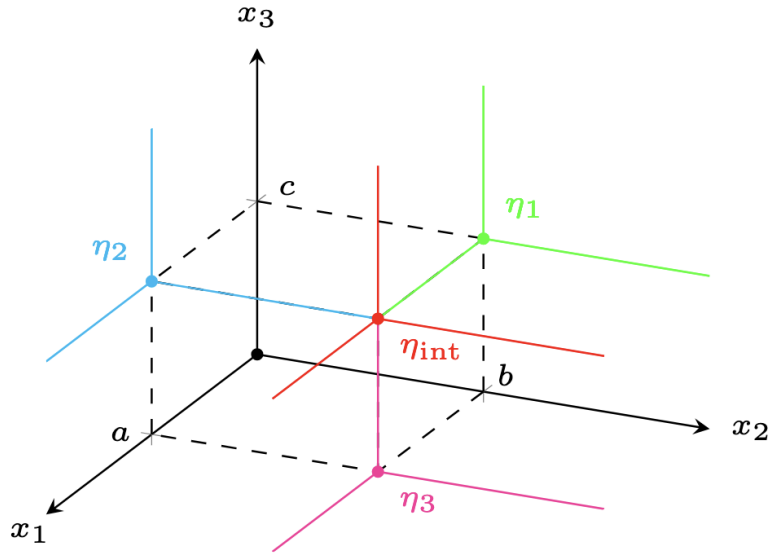}
\caption{Basepoints of plane partitions $\eta_1, \eta_2, \eta_3$ in $\R^3$}
\label{fig:pisinR3}
\end{figure}

Consider triples of plane partitions $\eta = (\eta_1, \eta_2, \eta_3)$ such that $\eta_1$ is based at $(0,b,c)$, $\eta_2$ is based at $(a,0,c)$, and $\eta_3$ is based at $(a,b,0)$ in $\R^3$ (see Figure \ref{fig:pisinR3}). We say that a box $(i,j,k)$ is in the \emph{intersection space} if $i \geq a$, $j \geq b$, and $k \geq c$. We denote all the boxes in the intersection space by $\eta_{\text{int}}$.

We define different types of boxes based on the number of plane partitions they are contained in as:

\begin{definition}
We say that a box $(i,j,k) \in \eta = (\eta_1, 
\eta_2, \eta_3)$ is:
\begin{itemize}
    \item a type III box if $(i,j,k) \in \eta_1, \eta_2, \eta_3$ (triple intersection boxes)
    \item a type II box if $(i,j,k) \in \eta_m, \eta_n$ and $(i,j,k) \notin \eta_l$ for $\{m,n,l\} = \{1,2,3\}$ (double intersection boxes)
    \item a type I box if $(i,j,k) \in \eta_{m}$ and $(i,j,k) \notin \eta_{n}, \eta_{l}$ for $\{m,n,l\} = \{1,2,3\}$ (boxes in only one of the plane partitions)
\end{itemize}
\end{definition}

Let $\eta_{\text{in}}$ denote the set of type III boxes, and let $\eta_{\text{out}}$ denote the set of type II boxes. Note that $\eta_{\text{in}} \cup \eta_{\text{out}}$ is plane partition based at $(a,b,c)$, with $\eta_{\text{in}} \cup \eta_{\text{out}} \subseteq \eta_{\text{int}}.$ We want to consider triples of plane partitions $\eta = (\eta_1, \eta_2, \eta_3)$ such that the following Criterion is satisfied:

\hspace{0.04cm}

\noindent \textbf{Criterion 1.} $\eta_{\text{int}} = \eta_{\text{in}} \cup \eta_{\text{out}}$,

\hspace{0.04cm}

\noindent that is, 

\begin{equation}
    \eta_{\text{int}} = (\eta_1 \cap \eta_2) \cup (\eta_1 \cap \eta_3) \cup (\eta_2 \cap \eta_3).
\end{equation}

We define an equivalence relation on triples of plane partitions satisfying Criterion 1 as follows. If $\eta = (\eta_1, \eta_2, \eta_3)$ and $\tilde{\eta} = (\tilde{\eta_1}, \tilde{\eta_2}, \tilde{\eta_3})$, we say that $\eta \sim \tilde{\eta}$ if they have the same multiset of boxes. That is, $\eta \sim \tilde{\eta}$ if:

\begin{itemize}
    \item $\eta_{\text{in}} = \tilde{\eta}_{\text{in}}$ (type III boxes the same)
    \item $\eta_{\text{out}} = \tilde{\eta}_{\text{out}}$ (type II boxes the same, regardless of which two partitions they came from)
    \item $\eta_1$ agrees with $\tilde{\eta}_1$ on $[0,a) \times [b, \infty) \times [c, \infty)$ 
    \item $\eta_2$ agrees with $\tilde{\eta}_2$ on $[a, \infty) \times [0, b) \times [c, \infty)$ 
    \item $\eta_3$ agrees with $\tilde{\eta}_3$ on $[a, \infty) \times [b, \infty) \times [0,c)$ 
\end{itemize}

\noindent The last three conditions ensure that all type I boxes, that is, those not in the intersection space by Criterion 1, are the same. We are now ready to define double-box configurations as:

\begin{definition}\label{def:double_box_config}
    Given $(a,b,c) \in \N^3$, an equivalence class of triples of plane partitions $\eta = (\eta_1, \eta_2, \eta_3)$ satisfying Criterion 1 under the equivalence relation $\sim$ is called a double-box configuration. Note that we often denote such an equivalence class by $\eta$, rather than $[\eta]$.
\end{definition}

We denote the set of all double-box configurations by $DBC_{a,b,c}$. For each $\eta \in DBC_{a,b,c}$, we define the following:

\begin{definition}\label{def:size_of_DBC}
The weight of a double-box configuration $\eta= (\eta_1, \eta_2, \eta_3)$ is defined as
$$|\eta| = |\eta_1| + |\eta_2| + |\eta_3| - |\eta_{\text{int}}|$$
$$ = \#\{\text{type $I$ boxes}\} + \#\{\text{type II boxes}\} + 2 \cdot \#\{\text{type III boxes}\}.$$
Note that this quantity is well-defined on equivalence classes.
\end{definition}

Next, we consider elements within the equivalence class of a double-box configuration $[\eta]$. To do this, we first make the following definition:

\begin{definition}
    A box $(i,j,k) \in \eta_{\text{out}}$ is said to be moveable if there exists $\hat{\eta} \neq \tilde{\eta} \in [\eta]$ and two indices $m \neq n \in \{1,2,3\}$ such that $(i,j,k) \notin \hat{\eta}_m$ and $(i,j,k) \notin \tilde{\eta}_n$. 
\end{definition}

\noindent Triples of plane partitions within an equivalence class (i.e. a double-box configuration), may differ by which two plane partitions a moveable box is contained in. We consider several examples to illustrate this.

\begin{ex}
    Let $(a,b,c) = (1,1,1)$. Consider $\eta = (\eta_1, \eta_2, \eta_3)$ and $\tilde{\eta} = (\tilde{\eta}_1, \tilde{\eta}_2, \tilde{\eta}_3)$ as defined in Figure \ref{fig: equiv_class_example}.

    \begin{figure}[H]
        \centering
        \includegraphics[width=\textwidth]{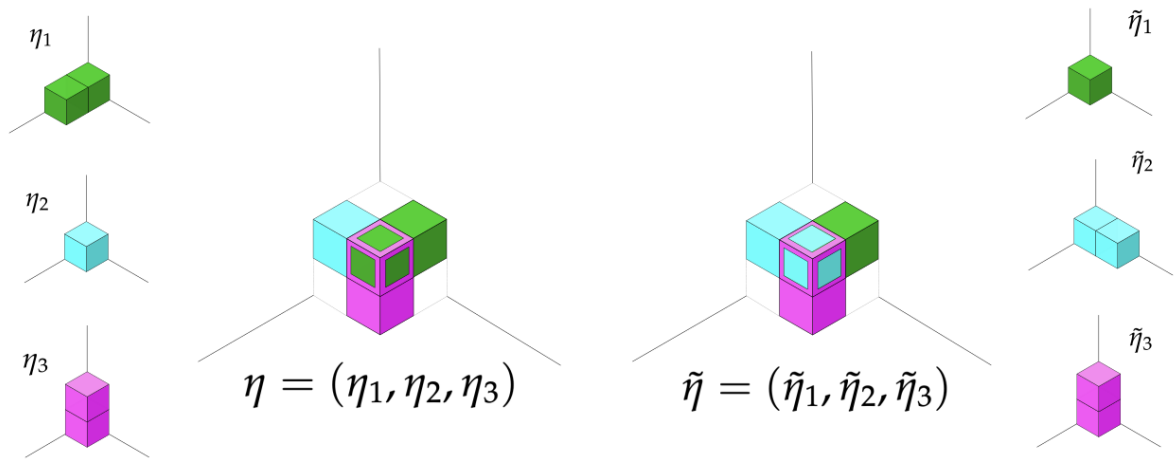}
        \captionof{figure}{Example of $\eta \neq \tilde{\eta}$ with $[\eta] = [\tilde{\eta}]$.}
        \label{fig: equiv_class_example}
    \end{figure}
    
    \noindent There is one type II box in this double-box configuration at $(1,1,1)$. In $\eta$, this box is contained in $\eta_1$ and $\eta_3$, and in $\tilde{\eta}$ this box is contained in $\tilde{\eta}_2$ and $\tilde{\eta}_3$. Since $\eta$ and $\tilde{\eta}$ contain the same multiset of boxes, we have that $[\eta] = [\tilde{\eta}]$, and so the type II box at $(1,1,1)$ is a moveable box.
\end{ex}

\newpage 

\begin{ex}
    Let $(a,b,c) = (1,1,1)$, and consider $\eta = (\eta_1, \eta_2, \eta_3)$ as defined in Figure \ref{fig: DBC_example}.

    \begin{figure}[H]
        \centering
        \includegraphics[width=0.5\textwidth]{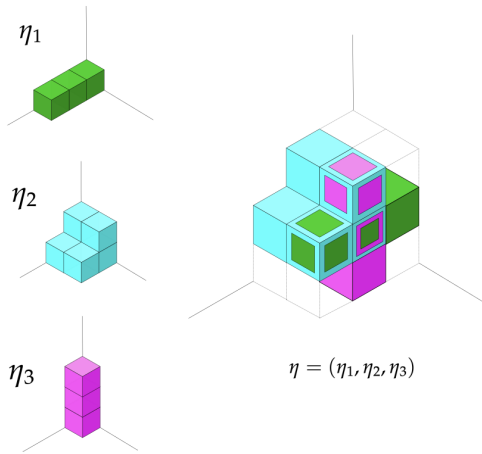}
        \captionof{figure}{Example of a double-box configuration.}
        \label{fig: DBC_example}
    \end{figure}

\noindent In this example, there is one type III box at $(1,1,1)$, and two type II boxes, one at $(2,1,1)$ contained in $\eta_1$ and $\eta_2$, and one at $(1,1,2)$ contained in $\eta_2$ and $\eta_3$. Both of these type II boxes are not moveable. The box at $(1,1,2)$ cannot be contained in $\eta_1$ because $(0,1,2) \notin \eta_1$, and the box at $(2,1,1)$ cannot be contained in $\eta_3$ because $(2,1,0) \notin \eta_3$. 
\end{ex}

To define the generating function for double-box configurations, we first make the following definition:

\begin{definition}\label{def:contribution_of_DBC}
The contribution of a double-box configuration $\eta$ is defined as

$$\chi(\eta) = 2^m$$

\noindent where $m$ is the number of connected components of moveable boxes in $\eta$. Two boxes are in the same connected component if they share a face.
\end{definition}

Finally, we define the generating function for double-box configurations as

\begin{equation}\label{eq:DBC_GF}
    Z_{DBC}^{a,b,c}(q) = \sum_{\eta \in DBC_{a,b,c}}\chi(\eta)q^{|\eta|}, 
\end{equation}

\noindent where $|\eta|$ is defined in Definition \ref{def:size_of_DBC}, and $\chi(\eta)$ is defined in Definition \ref{def:contribution_of_DBC}.

\subsection{Tripartite double-dimer configurations}

\begin{figure}[H]
    \centering
    \includegraphics[width=0.5\textwidth]{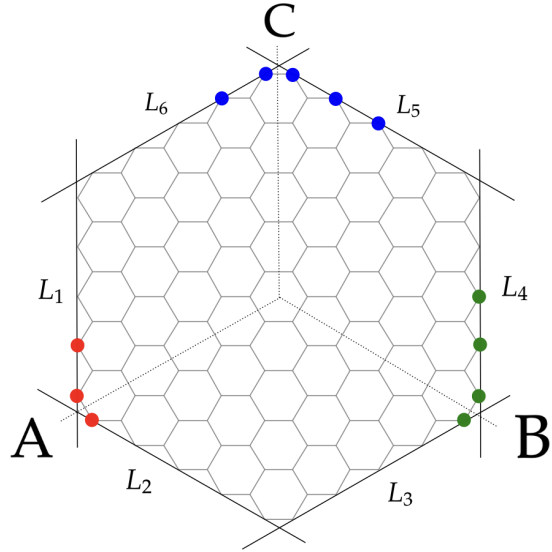}
    \captionof{figure}{Red, green, and blue nodes on $H(5)$, with $a = 2$, $b=3$, and $c=1$.}
    \label{fig:nodes}
\end{figure}

In this section we will define the tripartite node pairing of the double-dimer configurations on the hexagon graph. Let $H(n)$ be the hexagon graph of size $n \times n \times n$.  That is,  project the points $\{0, \dots, n\}^3 \subset \mathbb{N}^3$ onto the plane $P: \{x+y+z=0\}$ to obtain the vertices of a hexagon-shaped piece of the triangular lattice; $H(n)$ is the planar dual of this graph, without an external vertex.  Choose coordinates $(x, y)$ for P such that a third of the edges are parallel to the $x$ axis - "horizontal" - and the others have slope $\pm \sqrt{3}/2$.  For convenience, we will use standard "compass coordinates" to describe directions on this picture - so "North" is the positive $y$ direction, "West" is negative $x$, and so on  
 (see Figure~\ref{fig:nodes}).
 



Next, let $A$ be the southwest corner of $H(n)$ (that is, the intersection of the lines $L_1$ and $L_2$ in Figure \ref{fig:nodes}), let $B$ be the southeast corner (intersection of $L_3$ and $L_4$), and let $C$ be the north corner of $H(n)$  (intersection of $L_5$ and $L_6$). We define the following sets of nodes, i.e. special vertices, on the boundary of $H(n)$:

$$R = \{\text{$a$ nodes on $L_1$ closest to $A$}\} \cup \{\text{$c$ nodes on $L_2$ closest to $A$}\},$$

$$G = \{\text{$c$ nodes on $L_3$ closest to $B$}\} \cup \{\text{$b$ nodes on $L_4$ closest to $B$}\},$$

$$B = \{\text{$b$ nodes on $L_5$ closest to $C$}\} \cup \{\text{$a$ nodes on $L_6$ closest to $C$}\}.$$

Note that $|R| = a+c$, $|G| = b+c$, and $|B| = a+b$, satisfy the triangle inequality, and so there is a unique planar \emph{tripartite pairing} of the nodes, we call this pairing $\sigma_{a,b,c}$. The planar tripartite pairing $\sigma_{a,b,c}$ matches the $a$ red nodes on $L_1$ with the $a$ blue nodes on $L_6$, the $b$ blue nodes on $L_5$ with the $b$ green nodes on $L_4$, and the $c$ red nodes on $L_2$ with the $c$ green nodes on $L_3$, so that each node is paired with another node of a different color. 

Let $DDC_{n}(\sigma_{a,b,c})$ be the set of all double-dimer configurations on $H(n)$ with node set $N = R \cup G \cup B$ and tripartite pairing $\sigma_{a,b,c}$. We define the generating function for elements in $DDC_n(\sigma_{a,b,c})$ as

\begin{equation}
    Z_{DDC}^{n; a,b,c}(q) = \frac{1}{w(\pi_0)}\sum_{\pi \in DDC_n(\sigma_{a,b,c})}2^{\ell(\pi)}w(\pi)
\end{equation}

\noindent where $\ell(\pi)$ is the number of closed loops of $\pi$, and the configuration $\pi_0 \in DDC_n(\sigma_{a,b,c})$ has minimal weight. The weight of a double-dimer configuration $\pi$ is given by

\begin{equation*}
w(\pi) = \prod_{e \in \pi}w(e)
\end{equation*}

\noindent where $w : E \to \mathbb{Q}[q]$ for an indeterminate $q$ is a weighting on the edge set of $H(n)$ such that the generating function for plane partitions is the same as the one for dimer configurations on $H(n)$ up to a constant.

We now make the following definition:

\begin{definition}\label{def:DDC} Let $DDC(\sigma_{a,b,c})$ be the set of all double-dimer configurations on the infinite hexagon lattice such that for each $\pi \in DDC(\sigma_{a,b,c})$, there exists $N \in \N$ such that $\pi$ restricted to $H(N)$ has the tripartite node pairing $\sigma_{a,b,c}$. 
\end{definition}


\noindent For elements of $DDC(\sigma_{a,b,c})$, we define the generating function as

\begin{equation}\label{eq:DDC_GF}
    Z_{DDC}^{a,b,c}(q) := \lim_{n \to \infty}Z_{DDC}^{n;a,b,c}(q) = \frac{1}{w(\pi_0)}\sum_{\pi \in DDC(\sigma_{a,b,c})}2^{\ell(\pi)}w(\pi).
\end{equation}

\noindent Proof that $Z_{DDC}^{a,b,c}(q)$ is well-defined is omitted here due to space constraints.

\section{Mapping double-box to double-dimer}

Our first main result is the following:

\begin{thm}\label{thm:GFs_equal}
Let $a,b,c \in \N$. Then

\begin{equation}
    Z_{DBC}^{a,b,c}(q) = Z_{DDC}^{a,b,c}(q)
\end{equation}

\noindent where $Z_{DBC}^{a,b,c}(q)$ is the generating function for double-box configurations as defined in Equation \ref{eq:DBC_GF}, and $Z_{DDC}^{a,b,c}(q)$ is the generating function for double-dimer configurations on the hexagon lattice with tripartite node pairing $\sigma_{a,b,c}$ as defined in Equation \ref{eq:DDC_GF}.
\end{thm}

\begin{figure}[H]
\centering
\includegraphics[width=0.5\textwidth]{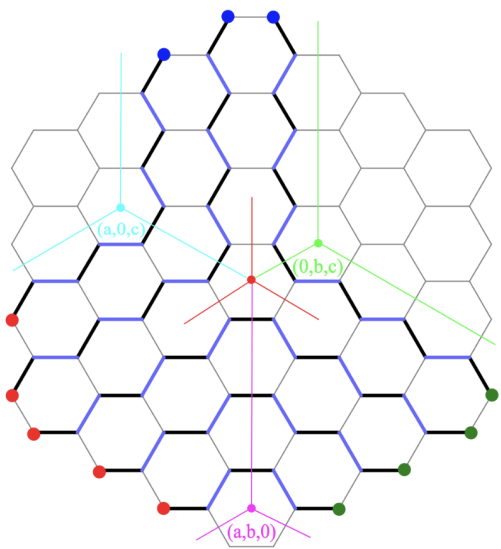}
\caption{Tripartite node pairing $\sigma_{a,b,c}$ with $a = 2$, $b = 1$, and $c = 3$. The red, green, and blue nodes on the outside face of the hexagon graph are connected such that no two nodes of the same color are endpoints of the same path. The projected basepoints via the folklore bijection of $\eta_1, \eta_2, \eta_3$ are in light green, light blue, and pink (compare to Figure \ref{fig:pisinR3}).}
\label{fig:nodes2}
\end{figure}

The mapping between double-box configurations and double-dimer configurations is easy to visualize via the folklore bijection, but significantly more difficult to prove. To visualize, we apply the folklore bijection (see Figure \ref{fig:folklore_bij}) to the double-box configurations. The result is a triple-dimer configuration on the hexagon graph, which is composed of a single-dimer configuration for each plane partition of $\eta = (\eta_1, \eta_2, \eta_3)$, overlayed such that the point $(a,b,c)$ in $\R^3$ is the center hexagon (see Figure \ref{fig:nodes2}). We show that removing the dimer configuration corresponding to the plane partition $\eta_{\text{int}}$ results in a double-dimer configuration on the hexagon lattice with the tripartite node pairing $\sigma_{a,b,c}$.

\section{Condensation relation}\label{sec:main_result}

Once we have the equivalence of generating functions given by Theorem \ref{thm:GFs_equal}, we may operate in the land of double-dimer configurations, which offers us many tools we can use to prove Theorem \ref{thm:main}. We show that both sides of Equation \ref{eq:main} satisfy a quadratic recurrence called \emph{condensation relation}, with the same initial conditions. 

On the right hand side of Equation \ref{eq:main}, we have $M(q)^2M_{a,b,c}(q)$, which we will denote by $X(a,b,c)$. We show that $X(a,b,c)$ satisfies the following recurrence relation

\begin{align*}
    X(a,b,c)X(a+1,b+1,c) &= X(a+1,b,c)X(a,b+1,c) \\
    & + q^{c}X(a+1,b+1,c-1)X(a,b,c+1).
\end{align*}

\noindent Note that $M(q)^4$ factors out of every term, and so we want to show that 

\begin{align*}
M_{a,b,c}(a)M_{a+1,b+1,c}(q) &= M_{a+1,b,c}(q)M_{a,b+1,c}(q) + q^cM_{a+1,b+1,c-1}(q)M_{a,b,c+1}(q).
\end{align*}

\noindent This follows from Kuo's graphical condensation (Theorem 6.2 in \cite{kuo2004applications}), with a slight modification (details omitted here), where $c$ is the height of the room.

For the left hand side of Equation \ref{eq:main}, we may consider $Z_{DDC}^{a,b,c}(q)$ instead of $Z_{DBC}^{a,b,c}(q)$ by Theorem \ref{thm:GFs_equal}. We apply \emph{double-dimer condensation}, a result of Jenne \cite{jenne2021combinatorics}, to $Z^{n;a,b,c}_{DDC}(q)$, which we denote here by $Z_{DDC}^n(a,b,c)$ to emphasize the change in $a,b,c$ in the recurrence. This gives the following relation

    \begin{align}
        Z_{DDC}^n(a,b,c)Z_{DDC}^n(a+1, b+1, c) & = Z_{DDC}^n(a, b+1, c)Z_{DDC}^n(a+1, b, c) \label{eq: trunc_DDC_recurrence}\\
         & +  q^KZ_{DDC}^{n,down}(a+1, b+1, c-1)Z_{DDC}^{n,up}(a,b,c+1). \nonumber
    \end{align}

\noindent We calculate this constant $K$, and show that it does not depend on $n$. Thus the limit as $n$ goes to infinity is well-defined, and so $Z_{DDC}^{a,b,c}(q)$ also satisfies this recurrence. Lastly, we show that both sides satisfy the same initial conditions. Details are omitted here, due to space constraints.

\acknowledgements{We would like to thank Amin Gholampour and Martijn Kool for helpful conversations, as well as Cruz Godar and Ava Bamforth.  Some illustrations in this paper were created with Dimerpaint, written by Leigh Foster and the second author.}

\bibliographystyle{alpha}
\bibliography{dd}

\end{document}